\input phyzzx
\def\e{\adveq\eqno{\rm (\chapterlabel\the\equanumber)}}

\def\adveq{\global\advance\equanumber by 1}
\def\myeq{{\rm \chapterlabel\the\equanumber}}

\def\d#1{{\rm d}#1\,}

\def\semidirect{\mathrel{\raise0.04cm\hbox{${\scriptscriptstyle |\!}$
\hskip-0.175cm}\times}}


\def\ref#1{$^{[#1]}$}

\def\r#1{$[\rm#1]$}

\def\d{{\rm d}}

\def\e{\adveq\eqno{\rm (\chapterlabel\the\equanumber)}}

\def\adveq{\global\advance\equanumber by 1}
\def\myeq{{\rm \chapterlabel\the\equanumber}}

\def\d#1{{\rm d}#1\,}

\def\semidirect{\mathrel{\raise0.04cm\hbox{${\scriptscriptstyle |\!}$
\hskip-0.175cm}\times}}


\def\ref#1{$^{[#1]}$}

\def\r#1{$[\rm#1]$}

\def\d{{\rm d}}

\def\conjecture#1{{\rm Conjecture: #1}}
\date{February, 2005}
\titlepage
\title{On The Distribution of Prime Multiplets (II)}
\line{\hfil Doron Gepner \hfil}
\line{\hfil Weizmann Institute of Science \hfil}
\line{\hfil Rehovot 76100, Israel \hfil}
\abstract
Recently, I have defined the so called PDF's (prime distribution factors)
 which govern the distribution
of prime numbers of the type $p,p+a_i$ being all primes up to some number
$n$. It was shown that the PDF's are expressible in terms of the basic PDF's 
which are defined as $a_i-a_j$ or $a_i$ being composed of primes which
are less or equal to the number of primes. For example, $p,p+2$ 
(twin primes), or $p,p+2,p+6$ being all primes (basic triplets).
We give here a conjecture for the number of basic prime PDF's in terms
of Hardy-Littlewood numbers, thus completing the determination of 
PDF's. These conjectures are supported by extensive calculations. 
\endpage
The distribution of primes has been fascinating since the work of
Gauss and Riemann. In ref. \REF\DG{Doron Gepner, On the distribution of
prime multiplets,
WIS preprint, May 2003, NT/0304477.}\r\DG\ the following observation
was made.
Define by $N(a_i,n)$ the number of sequences of $p,p+a_i, i=1\ldots m-1$
being all primes, up to the integer $n$. Then the following limit
exists and $C(a_i)$ is termed the "PDF" (prime distribution factor):
$$N(a_i,n)=C(a_i) \int_2^n {\d x \over \log(x)^m},\e$$
for large enough $n$. The PDF's were shown to be expressible
in terms of the basic PDF's by a product with some rational number.
For example, the PDF for $p$ and $p+n$ ($n$ even) is given by
$$C(p,p+n)=C(p,p+2) \prod_{q>2 \ {\rm prime} \atop q|n} {q-1\over q-2}.\e$$

The basic PDF's are defined by the requirement that $a_i-a_j$
and $a_i$ are all products of primes less or equal to $m$.
For example: $p,p+2$, for twins, $p,p+2,p+6$, for triplets,
$p,p+2,p+6,p+8$ for quadruplets, $p,p+2,p+6,p+8,p+12$ for
quintuplets, add $p+18$
for sextuplets, etc.

A conjecture for prime twins, $p,p+2$ is known before, see, e.g.,
\REF\Br{R.P. Brent, Math. Comp. 28 (1974)}\r\Br. It is given by
$$C(2)=2 c_2,\e$$
where
$$c_2=\prod_{p>2 \atop p\ {\rm prime}}^\infty 1-{1\over (p-1)^2} ,\e$$
and is numerically $C(2)=1.32032...$

Our generalization to this is to use Hardy-Littlewood numbers,
defined by
$$c_m=\prod_{p>m \atop p\ {\rm prime}}^\infty {p^{m-1} (p-m)\over (p-1)^m},\e $$
The first few HL (Hardy--Littlewood) numbers are: $c_2=0.6601618$,
$c_3=0.6351663...$, $c_4=0.3074948...$, $c_5=0.409874...$,
etc. See, e.g., ref.
\REF\HT{P. Moree et. al., http://www.gn-50uma.de/alula/essays/Moree/Moree.en.shtml}
\r\HT, where these numbers are calculated, and references therein.
Our main observation is that the basic PDF's are given by the 
HL numbers times a "simple" rational number:

\conjecture
The basic PDF's are given by HL numbers times a rational number.
This number is $2$ for twins, $9/2$ for triplets, $27/2$ for
quadruplets.

Thus, we have for the basic PDF's the formula:

$$C(2)=2 c_2=1.32032...\e$$
$$C(3)=9 c_3/2=2.858248...\e$$
$$C(4)=27 c_4/2=4.1511808...\e$$

To check these conjectures we made an extensive calculation of the
number of basic $m$-plets. The program we used was written in MATLAB
and is listed in the appendix.

For twin primes we find up to $10^{10}$, $27,412,673$ twins, $p$ and
$p+2$ being primes. This gives the PDF $1.32038...$ which is
$4.5\times 10^{-5}$ off the limit $C(2)$, eq. (6).
For up to $n=2\times 10^{10}$ we find $4,942,554$ basic triplets, of the
form $p,p+2,p+6$, giving the
PDF $2.85768$ which is $2 \times 10^{-4}$ off $C(3)$, eq. (7).
For basic quadruplets, $p,p+2,p+6,p+8$, we find $898,998$ up to the number
$7\times 10^{10}$,
which gives the PDF, $4.1503..$, which is $2 \times 10^{-4}$ off 
$C(4)$, eq. (8). This supports the above conjecture.

For basic quintuplets we were unable  to determine the rational number,
which seems not as simple as for $m\leq 4$. For $n=4 \times 10^{11}$ we find
$370,502$ basic quintuplets, which gives the PDF, $10.1193$ or
$C(5)\approx 24.6888 \times c_5$. Thus, we know that the rational factor
is $24.6888$ approximately, but it is not enough to form a conjecture
about it.

To summarize, we have a conjecture about the number of all prime multiplets,
expressing them in terms of Hardy--Littlewood numbers, times some
known rational coefficients. These completes the calculation
of the number of prime multiplets up to quadruplets.
For higher multiplets, we were unable to determine the coefficients,
which seems not so simple, and it is left to further work.

\refout

\appendix

Attached are two MATLAB functions used to calculate the number of
basic prime multiplets.
We used the function, "primel2", to determine the number of primes
between n1 and n2, employing the sieve method. Further below is a sample 
function, "pairs5",  which 
calculates the number of basic quintuplets up to the integer STEP$\times$ mm,
which works by looking on consequent primes from a list generated
by "primel2", and determining
their differences.

\obeylines

  function p = primel2(n1,n2)
\%PRIMEL2 Generate list of prime numbers.
\%   PRIMEL2(N1,N2) is a row vector of the prime numbers less than or 
\%   equal to N2 and bigger or equal to N1.   
p1=primes(sqrt(n2));
lp1=length(p1);
qq=n1:n2;
if n1==1; qq(1)=0;
end;
for mm=1:lp1
    qw=p1(mm);
    sw2=n1-mod(n1-1,qw)+qw-1;
    if sw2==qw ;sw2=sw2+qw;
    end;
        qq(sw2-n1+1:qw:n2-n1+1)=0;
end;
p=qq($qq>0$);
\%
\%
\% Below is the sample function, "pairs5", which calculates the number
\% of basic quintuplets, up to the integer $mm\times step$. The reason of doing
\% this in "batches" is to prevent exceeding the available memory.
\% We typically take step=$10^7$.
\%
function np =  pairs5(mm,step)
np=0;
for jj=1:mm;
    n1=step*(jj-1)+1;
    n2=step*jj;
    if $n1>1$; n1=n1-11;
    end;
    pp=primel2(n1,n2);
   np=np+sum([0 0 0 pp 0]-[0 0 0 0 pp]==2 \& [0 0 pp 0 0]-[0 0 0 pp 0]==4 ...
   \& [0 pp 0 0 0]-[0 0 pp 0 0]==2 \& [pp 0 0 0 0]-[0 pp 0 0 0]==4);
end;        

\bye